\documentclass[12pt]{amsart}

\author{Stefano Decio}
\address{Department of Mathematical Sciences, Norwegian University of Science and Technology, 7491 Trondheim, Norway}
\email{stefano.decio@ntnu.no}

\usepackage{amsmath}
\usepackage{amssymb}
\usepackage{amsfonts}
\usepackage{enumerate}
\usepackage{mathrsfs}
\usepackage{MnSymbol}
\usepackage{mathtools}

\newcommand{\R}{{\mathbf R}}

\newcommand{\ld}{{\lambda}}

\newtheorem{theorem}{Theorem}
\newtheorem{lemma}{Lemma}
\newtheorem{claim}{Claim}

\newtheorem{proposition}{Proposition}
\theoremstyle{remark}
\newtheorem*{remark}{Remark}
\theoremstyle{definition}

\begin{document}

\begin{abstract}
    We show that Steklov eigenfunctions in a bounded Lipschitz domain have wavelength dense nodal sets near the boundary, in contrast to what can happen deep inside the domain. Conversely, in a two-dimensional Lipschitz domain $\Omega$, we prove that any nodal domain of a Steklov eigenfunction contains a half-ball centered at $\partial\Omega$ of radius $c_{\Omega}/\ld$.
\end{abstract}

\title{Nodal sets of Steklov eigenfunctions near the boundary: Inner radius estimates}
\title[Inner radius estimates]{Nodal sets of Steklov eigenfunctions near the boundary: Inner radius estimates}
\maketitle

\section{Introduction}
The Dirichlet-to-Neumann operator, which maps a function on the boundary of a domain to the normal derivative of its harmonic extension in the domain, is an object with a revered history in spectral theory, its eigenvalues and eigenfunctions still being the subject of active research. In this note we study the harmonic extensions of eigenfunctions of the Dirichlet-to-Neumann problem, usually called Steklov eigenfunctions; these have attracted a good amount of attention in the past years in their own right, as an interior counterpart to the eigenfunctions of the Dirichlet-to-Neumann operator. 

For $\Omega$ a bounded domain in $\R^d$, where $d\geq 2$, the Steklov eigenvalue problem on $\Omega$ is informally given by
\begin{align}
\label{problem}
    \begin{cases}
    \Delta u=0 \qquad \: \text{in} \ \Omega,\\
    \partial_{\nu}u=\ld u \quad \ \text{on} \ \partial \Omega. \end{cases}
\end{align}
Note that the spectral parameter $\lambda$ is in the boundary conditions; a number $\ld$ for which there is a solution to \eqref{problem} is called a Steklov eigenvalue, and a corresponding solution is called a Steklov eigenfunction. Of course, for problem \eqref{problem} to make sense some regularity on $\partial\Omega$ must be imposed. In most of this work we will require that $\Omega$ be a Lipschitz domain; we defer the precise weak definition of problem \eqref{problem} in this setting to Section 2.\\

We refer to the survey \cite{GP} for many results and open questions about the Steklov problem, with an angle towards spectral geometry. In particular, a question was raised there about density of nodal (i.e. zero) sets of solutions of \eqref{problem}. It is a well known property that eigenfunctions of the Laplace-Beltrami operator on a compact Riemannian manifold have wavelength dense nodal sets, meaning that every ball of wavelength radius in the manifold contains a zero of the eigenfunction; the authors of \cite{GP} ask whether a similar property holds for the interior nodal sets of solutions of \eqref{problem}, or for nodal sets on the boundary, which are the nodal sets of eigenfunctions of the Dirichlet-to-Neumann operator. Note that the right density scale should be $C/\ld$; it is easy to convince oneself of this either by the scaling of \eqref{problem} or by looking at the explicit example of a ball.

A negative answer to the first question about interior nodal sets was given recently in \cite{BG}: the authors construct a dense family of simply connected domains in $\R^2$ with real analytic boundary such that for any such domain each Steklov eigenfunction $u_{\ld}$ does not change sign in a ball inside the domain, the location of the ball depending on $\ld$ but the radius being $\ld$-independent. They actually show more, namely that $u_{\ld}$ is well approximated in a fixed ball inside the domain by a function with only a finite number, independent of $\ld$, of non-zero Fourier modes. \\

The density question for boundary nodal sets is, to the best of our knowledge, in general still open; it follows from the main result of \cite{S}, Theorem 2.1, that it is true for smooth simply-connected planar domains. A closer inspection of the proof of Theorem 2.1 in \cite{S} reveals that $\mathcal{C}^{\infty}$-smoothness is not needed, and $\mathcal{C}^4$ should be sufficient. Note also that the boundary density property cannot possibly hold for every Lipschitz domain, or even for every piecewise smooth one, as the explicit example of a rectangle shows: eigenfunctions with zeros only on two sides of the rectangle exist (it is easy to compute the eigenfunctions explicitly as products of trigonometric and hyperbolic functions, see \cite{GP}). Our first contribution in the present note lies in between the interior and boundary cases: we show that balls of radius $C/\ld$ centered at the boundary do contain zeros of $u_{\ld}$. The precise statement follows.

\begin{theorem}
\label{upper}
Let $\Omega$ be a Lipschitz domain in $\R^d$, $d\geq 2$, and let $u=u_{\ld}$ be a solution of \eqref{problem}, where we assume $\ld\neq 0$. Let $\mathcal{Z}_{u_{\ld}}=\{x\in \Omega : u_{\ld}(x)=0\}$. Then there exists a constant $C=C(\Omega)$ such that 
\begin{align}
\label{density}
    \mathcal{Z}_{u_{\ld}}\cap B\neq \emptyset
\end{align}
for any ball $B$ in $\R^d$ of radius $C/\ld$ centered at a point in $\partial\Omega$. 
\end{theorem}
Some remarks are in order. 
\begin{remark}
     \textit{(i)} The result is only interesting when $\ld$ is large enough: by orthogonality to constants on the boundary (that is, eigenfunctions with $\ld=0$) any solution of \eqref{problem} with $\ld\neq 0$ has to change sign on the boundary, and hence satisfy \eqref{density} for $C$ large enough if $\ld$ is bounded.\\
     \textit{(ii)} Our theorem does not exclude the chance that the zero set is essentially contained in a corona around the boundary of size shrinking with $\ld$; the results of \cite{BG} leave open the possibility that zeros are contained in any ball of radius $C/\ld$ centered at a point in a collar neighbourhood of the boundary of size independent of $\ld$, albeit small. We are inclined to think that the latter statement should be true.\\
     \textit{(iv)} The result of \cite{BG} says essentially that Steklov eigenfunctions need not oscillate deep in the interior, where they are small; there are also known instances where eigenfunctions concentrate on parts of the boundary, see for example \cite{LPPS}. Our result says that Steklov eigenfunctions need to oscillate near the boundary just because they satisfy the Steklov condition, regardless of whether they are small or large there. We thank Iosif Polterovich for bringing up this point. 

\end{remark}

\subsection*{Inner radius} 
We denote by $\mathcal{D}_{\ld}$ a nodal domain for the Steklov eigenfunction $u_{\ld}$, that is a connected component of the set $\Omega\setminus \mathcal{Z}_{u_{\ld}}$. Denote by $\rho(\mathcal{D}_{\ld})$ the radius of the largest ball $B$ centered at a point $y\in\partial \Omega\cap\overline{\mathcal{D}_{\ld}}$ such that $B\cap\Omega\subset \mathcal{D}_{\ld}$; we call it the modified inner radius of $\mathcal{D}_{\ld}$. Theorem \ref{upper} can then be restated as
\begin{align*}
    \rho(\mathcal{D}_{\ld})\leq \frac{C(\Omega)}{\ld}.
\end{align*}
Note that we do not say anything about balls in $\mathcal{D}_{\ld}$ centered at points not in $\partial\Omega$; as shown in \cite{BG}, it is possible that deep inside $\mathcal{D}_{\ld}$ one could inscribe balls of $\ld$-independent radius. We think that centering the balls at $\partial\Omega$ provides the correct analogy with the Laplace eigenfunctions case. In this setting, the second main contribution of this note is a sharp lower bound for the modified inner radius in dimension two.

\begin{theorem}
\label{lower}
Let $\Omega\subset \R^2$ be a Lipschitz domain, and let $\mathcal{D}_{\ld}$ be a nodal domain of $u_{\ld}$, where $\ld\neq 0$. There exists a constant $c(\Omega)>0$ such that 
\begin{align}
\label{inrad}
    \rho(\mathcal{D}_{\ld})\geq \frac{c(\Omega)}{\ld}.
\end{align}
\end{theorem}

In words, one can inscribe inside $\mathcal{D}_{\ld}$ the intersection of $\Omega$ and a ball of radius $c(\Omega)/\ld$ centered somewhere on $\partial \Omega\cap\overline{\mathcal{D}_{\ld}}$. Note that in this way we also get a lower bound for the ``genuine'' inner radius, by inscribing a full ball of radius $c_1/\ld$ inside the intersection of $\Omega$ and the ball above. 

\begin{remark}
     \textit{(i)} It is possible to obtain a non-sharp estimate (polynomial in $\ld^{-1}$) of the same type in higher dimensions as well, at least under some higher regularity assumption on the boundary ($\mathcal{C}^2$ suffices), by following the method of \cite{M}: this entails proving a local asymmetry estimate for Steklov eigenfunctions and combining it with a Poincaré inequality. As the flavour of the asymmetry estimate is somewhat different from the rest of the results in this paper, we do not pursue this here and defer it to later work. We will comment further on the issue in Section 4. 
\end{remark}

\subsection*{Hausdorff measure bounds}
In the Laplace eigenfunction case, a density result for zeros is the starting point in the proof of a sharp lower bound on the Hausdorff measure of the nodal set, conjectured by Yau. The result is due to Logunov in \cite{L}: the proof is not easy and uses some clever insights about harmonic functions. In upcoming work, we combine Theorem \ref{upper} with the ideas of \cite{L} to improve on the existing lower bound for the Hausdorff measure of the nodal set of a Steklov eigenfunction: namely, we prove that if $\Omega$ is a domain with $\mathcal{C}^2$-smooth boundary, for a Steklov eigenfunction $u_{\ld}$ the bound 
\begin{align*}
    \mathcal{H}^{d-1}(\mathcal{Z}_{u_{\ld}})\geq c(\Omega)
\end{align*}
holds with a constant $c(\Omega)$ independent of $\ld$. The previous best result, established in \cite{SWZ}, yielded a lower bound polynomially decaying in $\ld$ (for $d>2$), of order $\ld^\frac{2-d}{2}$. In the same article we also obtain an upper bound on the Hausdorff measure of the nodal set. 

\subsection*{Plan of the paper} 
We prove Theorem \ref{upper} in Sections 2 and 3; in Section 2 we discuss some regularity for Steklov eigenfunctions and we describe an auxiliary boundary value problem that we need in the argument, and Section 3 gives the actual proof of the density statement (two of them, as a matter of fact). The proof of Theorem \ref{lower} is given in Section 4.

\section{On boundary value problems}
\subsection*{The Steklov problem}
In the statement of Theorem \ref{upper} we assume that $\Omega$ is a Lipschitz domain: by this we mean that for every point $p\in \partial \Omega$ there is a hyperplane $H$ through $p$ such that locally $\partial \Omega$ is the graph of a Lipschitz function on $H$ and $\Omega$ is the set of points below that graph. In this setting, a weak solution of \eqref{problem} is a function $u\in H^1(\Omega)$ such that 
\begin{align}
\label{weak}
    \int_{\Omega}\nabla u \nabla v = \ld \int_{\partial\Omega} uv
\end{align}
for all $v\in H^1(\Omega)$ (we denote by $H^1(\Omega)$ the Sobolev space $W^{1,2}(\Omega)$); a number $\ld$ for which the above equation has a nontrivial solution is called an eigenvalue of the Steklov problem, and a solution $u=u_{\ld}$ is a Steklov eigenfunction. Since $\Omega$ is Lipschitz, the trace operator $H^1(\Omega)\rightarrow L^2(\partial\Omega)$ is compact and hence the spectrum of \eqref{problem} is discrete and accumulating to infinity; the eigenvalues can be arranged in a sequence $0=\ld_0<\ld_1\leq \ld_2\leq\dots\rightarrow \infty$, with corresponding eigenfunctions $u_{\ld_j}$ forming an orthonormal basis of $L^2(\partial\Omega)$. This is a folklore result, but see \cite{AM} for a proof. 

It is possible to prove, using a Moser iteration type argument, that a solution of \eqref{weak} is actually in $L^{\infty}(\Omega)$, with bounds for the $\sup$ norm depending on $\ld$ and on $\Omega$ (see for instance \cite{BGT}). One thus sees that a Steklov eigenfunction is the solution of a Neumann boundary value problem with Neumann data in $L^{\infty}(\partial\Omega)$; in \cite{BH}, a solution to the Neumann problem with bounded Neumann data is constructed, and this solution is continuous up to the boundary. By the uniqueness up to constants of weak solutions of Neumann boundary value problems in Lipschitz domains (see \cite{JK}), one obtains then that solutions to \eqref{weak} actually are in $\mathcal{C}(\overline{\Omega})$. This fact is not really needed in the next section, but in our opinion it makes it easier to think about Steklov eigenfunctions in Lipschitz domains.

\subsection*{A Steklov-Dirichlet problem}
Our strategy to prove density for the zero sets of Steklov eigenfunctions near the boundary relies on an auxiliary boundary value problem, for which we need a description of the spectrum analogous to that of the Steklov problem. 

Let $W$ be a bounded Lipschitz domain in $\R^d$; we assume the boundary consists of two non-empty relatively open sets, $F$ and $S$, with common Lipschitz boundary. We call a Steklov-Dirichlet (sometimes called mixed Poincar{\'e}-Steklov as in \cite{A}) problem the following mixed spectral problem:
\begin{align}
\label{mixed}
    \begin{cases}
    \Delta u=0 \qquad \: \text{in} \ W,\\
    u=0 \qquad \ \ \; \text{on} \ S,\\
    \partial_{\nu}u=\sigma u \quad \ \text{on} \ F.\\
    \end{cases}
\end{align}
Under our hypotheses, it is well known that the spectrum is discrete and consists of eigenvalues $0<\sigma_1\leq \sigma_2\leq\dots \rightarrow \infty$, each having finite multiplicity (see \cite{A}, \cite{BKPS}); moreover, the corresponding eigenfunctions are orthogonal in $L^2(F)$. For completeness of the exposition we provide an outline of a proof, following the arguments of \cite{AM} and \cite{AM2} for the Steklov case. As in \cite{A}, \cite{BKPS}, we denote
$$H_0^1(W;S)=\{u\in H^1(W) : u=0 \ \text{on} \ S\};$$
here and in the sequel, we consider elements in $H^1$ restricted to the boundary in the sense of Sobolev traces. We introduce now an operator of the Dirichlet-to-Neumann type: let $\mathcal{D}_{M}:L^2(F)\to L^2(F)$ be the unbounded operator defined by
\begin{align*}
    &dom(\mathcal{D}_{M})=\{\phi\in L^2(F) \ \text{such that} \ \exists u \in H_0^1(W;S) : \Delta u=0 \ \text{in} \ W ,\\ & \qquad \qquad \qquad u_{|_F}=\phi \ \text{and} \ \partial_{\nu}u\in L^2(F) \};\\
    &\mathcal{D}_{M}\phi=\partial_{\nu}u_{|_F}.
\end{align*}
Note that the spectrum of $\mathcal{D}_M$ coincides with that of problem \eqref{mixed}, so by general spectral theory one needs only to prove that $\mathcal{D}_M$ is selfadjoint, bounded below, and with compact resolvent. If $S=\emptyset$, in which case $\mathcal{D}_M$ is the standard Dirichlet-to-Neumann operator, this was proved in \cite{AM2},\cite{AM}. We follow closely their outline to check that the same holds when $S$ is non-empty. Let $V=\{u_{|_F} : u\in H_0^1(W;S)\}$ be the trace space, with $\|u_{|_F}\|_{V}=\|u\|_{H_0^1(W;S)}$; since $W$ is a Lipschitz domain the embedding $V\to L^2(F)$ is compact. We denote also 
\begin{align*}
    \mathcal{H}=\{u \in H_0^1(W;S) : \Delta u=0\}. 
\end{align*}
Since $0$ is not a Dirichlet eigenvalue, from the solvability of the Dirichlet problem in Lipschitz domains (see \cite{D},\cite{JK}) it follows that we have the (non-orthogonal) direct sum 
$$H_0^1(W;S)=H_0^1(W)\oplus\mathcal{H}.$$
Using this, one sees that the linear map $\mathcal{H}\ni u\to u_{|_F}\in V$ is bijective. Define now the bilinear form $a : V\times V\to \R$ by 
$$a(\phi,\psi)=\int_{W}\nabla u \nabla v,$$
where $u,v \in \mathcal{H}$ with $u_{|_F}=\phi$ , $v_{|_F}=\psi$. We recall now an abstract lemma, which can be found as Lemma 2.3 in \cite{AM2}:
\begin{lemma}
    Let $X_1,X_2,X_3$ be Banach spaces, $X_1$ reflexive. Let $T: X_1\to X_3$ be a compact linear operator and $S: X_1\to X_2$ a bounded injective linear operator. Then for any $\varepsilon>0$ there exists $c>0$ such that for any $x\in X_1$, $\|Tx\|^2_{X_3}\leq \varepsilon \|x\|^2_{X_1}+c\|Sx\|^2_{X_2}$.
\end{lemma}
Using the lemma for $T$ being the embedding $\mathcal{H}\to L^2(W)$ and $S$ being the trace operator $\mathcal{H}\to L^2(F)$, it is immediate to show that $a$ is coercive (elliptic) on $V$. There is a linear operator on $L^2(F)$ associated to $a$, call it $B$; let $u\in \mathcal{H}$, $b\in L^2(F)$. Then $u_{|_F}\in dom(B)$ and $Bu_{|_F}=b$ if and only if 
\begin{equation}
\label{form}
    \int_{W}\nabla u \nabla v=\int_{F} bv
\end{equation}
for all $v\in\mathcal{H}$. Using that $\Delta u=0$, one gets that \eqref{form} holds also for $v\in H_0^1(W)$; remembering that $H_0^1(W;S)=H_0^1(W)\oplus\mathcal{H}$, \eqref{form} holds for all $v\in H_0^1(W;S)$, and consequently $b=\partial_{\nu} u$; this implies that $u_{|_F}\in dom(\mathcal{D}_M)$ and $\mathcal{D}_M u_{|_F}=Bu_{|_F}$. It is also easy to see that $dom(\mathcal{D}_M)\subset dom(B)$, so that $\mathcal{D}_M=B$ is the operator associated with the form $a$. The fact that $\mathcal{D}_M$ is bounded below and with compact resolvent now follows immediately by coercivity of $a$ and the embedding $V\to L^2(F)$ being compact. \\

The description of the spectrum of $\mathcal{D}_M$ is now given by the usual max-min procedure with the form $a$; in particular, the first eigenvalue of $\mathcal{D}_M$, and thus of problem \eqref{mixed}, is given by
\begin{align}
\label{eigenvalue}
    \sigma_1=\min_{{u\in H_0^1(W ; S)}}\frac{\int_W |\nabla u|^2}{\int_F |u|^2}.
\end{align}
Let now $u_j$, $u_k$ correspond to eigenvalues $\sigma_j\neq \sigma_k$; one has the string of equalities
\begin{align*}
    \sigma_j\int_F u_j u_k = \int_F \partial_{\nu}u_j u_k = \int_W \nabla u_j \nabla u_k = \int_{F} u_j \partial_{\nu}u_k = \sigma_k\int_F u_j u_k,
\end{align*}
which implies that $u_j$ and $u_k$ are orthogonal in $L^2(F)$.
We record the following fact, which is proved in the same way as in the Dirichlet case:
\begin{proposition}
    The first eigenvalue $\sigma_1$ of problem \eqref{mixed} is simple and the corresponding eigenfunction is positive in $W$.
    \begin{proof}
    Let $u$ be an eigenfunction corresponding to $\sigma_1$. Since $u$ is in $H_0^1(W;S)$, also $|u|\in H_0^1(W;S)$. Moreover, $|\nabla|u||=|\nabla u|$, so that the variational quotient in \eqref{eigenvalue} is minimized by $|u|$ as well; as a consequence, $|u|$ is harmonic in $W$, and therefore $u$ does not change sign in $W$. Without loss of generality, we can assume $u$ is positive. If there were two linearly independent eigenfunctions $u_1$ and $u_2$ corresponding to $\sigma_1$, then we could find a linear combination of them which changes sign in $W$. By the above, this implies that $\sigma_1$ is simple.
    \end{proof}
\end{proposition}

\section{Density of nodal sets near the boundary}
We give here two different proofs of Theorem \ref{upper}. Throughout the section we denote by $u=u_{\ld}$ a solution of problem \eqref{problem}, corresponding to the eigenvalue $\ld$. We can assume that $\ld$ is greater than some large enough constant $\ld_{\Omega}$ (see Remark \textit{(i)} following Theorem \ref{upper}). The two proofs are quite similar in spirit, and both involve the mixed Steklov-Dirichlet problem introduced in Section 2; while the second proof may look more natural to those familiar with Laplace eigenfunctions, it appears to require a higher degree of regularity on the domain. The first proof works under the Lipschitz assumption only and, moreover, one does not really need to think about nodal domains for the argument to carry on.
\subsection*{First proof} Let $x_0$ be a point on $\partial\Omega$, and let $W=B(x_0,s)\cap \Omega$. We take $s$ small enough so that $W$ can be written as the set of points below the graph of a Lipschitz function. Assume that $u$ does not change sign in $W$, say $u$ is positive. To keep the notation of Section 2, we call $S=\partial B(x_0,s)\cap \Omega$ and $F=\partial\Omega\cap B(x_0,s)$. We let $\sigma_1$ and $h$ denote respectively the first eigenvalue and the first eigenfunction of the Steklov-Dirichlet problem in $W$ with boundaries $S$ and $F$; that is, $h$ solves the problem
\begin{align}
\label{auxiliary}
    \begin{cases}
    \Delta h=0 \qquad \ \: \text{in} \ W,\\
    h=0 \qquad \ \ \ \  \text{on} \ S,\\
    \partial_{\nu}h=\sigma_1 h \quad \ \text{on} \ F.\\
    \end{cases}
\end{align}
The key to the proof is the following lemma, which is a consequence of the variational principle \eqref{eigenvalue}. We warn the reader that constants denoted by $C$ may change value from line to line.
\begin{lemma}
\label{key}
There exists a small enough constant $c_{\Omega}$ and a dimensional constant $C>0$ such that if $s<c_{\Omega}$, then $\sigma_1<C/s$.
\begin{proof}
We select $c_{\Omega}$ such that for $s<c_{\Omega}$,
\begin{align}
\label{measure}
    \mathcal{H}^{d-1}(F\cap B(x_0,s/2))\gtrsim s^{d-1},
\end{align}
where $\mathcal{H}^{d-1}$ denotes $d-1$-dimensional Hausdorff measure and the symbol $\gtrsim$ hides a dimensional constant. Here we used the hypothesis that the domain is Lipschitz. Consider now the function
\begin{align*}
    \psi(x)=1-\frac{|x-x_0|^2}{s^2},
\end{align*}
and note that $\psi\in H_0^1(W ; S)$ so that $\psi$ is allowed as a test function in the variational quotient in \eqref{eigenvalue}. We have that $$\int_W |\nabla \psi|^2\lesssim \frac{1}{s^4}\int_{B(x_0,s)}|x-x_0|^2\lesssim s^{d-2}.$$
Note that for $x\in W\cap B(x_0,s/2)$ we have $\psi(x)\geq 3/4$, so that by \eqref{measure} 
$$\int_{F}|\psi|^2\geq \int_{F\cap B(x_0,s/2)}|\psi|^2\gtrsim s^{d-1}.$$
Plugging the above estimates in \eqref{eigenvalue}, we obtain $\sigma_1\lesssim 1/s$.
\end{proof}
\end{lemma}

We take now the function $v=h^2$, where $h$ is as in \eqref{auxiliary}. Note that $\Delta v=2|\nabla h|^2\geq 0$, $v=\partial_{\nu}v=0$ on $S$ and $\partial_{\nu}v=2\sigma_1 v$ on $F$. By Green's formula (see for instance \cite{MZ}, Chapter 6, for a proof of Green's formula for Lipschitz domains) we thus have
\begin{align*}
    0=\int_{W}v\Delta u = \int_{W}u\Delta v + \int_F v\partial_{\nu}u - u\partial_{\nu} v.
\end{align*}
Since $u$ is positive in $W$, the first term on the right hand side is positive, so that by Lemma \ref{key} we obtain
\begin{align*}
    0\geq \int_F (\ld-2\sigma_1)uv \geq \int_F (\ld-C/s)uv,
\end{align*}
whence the condition $s<C/\ld$. The proof of Theorem \ref{upper} is complete. 

\subsection*{Second proof}
The main purpose in giving another proof is to record the observation contained in Proposition \ref{nodal}, which is certainly part of the folklore and has been used implicitly for instance in \cite{BKPS}. We assume here that $\partial\Omega$ has some higher regularity than Lipschitz, say $\partial\Omega$ is $\mathcal{C}^2$; under this hypothesis, a Steklov eigenfunction is in $\mathcal{C}^1(\overline{\Omega})$. Since we were unable to locate this fact in the literature, we state it as a proposition and outline a proof. 
\begin{proposition}
Let $\Omega$ be a $\mathcal{C}^{1,1}$ domain, and let $u=u_{\ld}$ be a Steklov eigenfunction. Then $u\in\mathcal{C}^1(\overline{\Omega})$.
\begin{proof}
    The argument relies on elliptic regularity results contained in Chapter 2 of \cite{Gri}; these results are used in \cite{LE} to prove analogous regularity statements for Neumann and Robin eigenfunctions. In particular, we will use the following statement, a special case of Proposition 2.5.2.3 in \cite{Gri}:\\
    \textit{Let $\Omega$ be a $\mathcal{C}^{1,1}$ domain, and let $u\in W^{1,p}(\Omega)$ be a weak solution of 
    \begin{align*}
    \begin{cases}
    -\Delta u+u=f \quad  \text{in} \ \Omega,\\
    \partial_{\nu}u=g \qquad \ \text{on} \ \partial \Omega, \end{cases}
    \end{align*}
    with $f\in L^p(\Omega)$ and $g\in W^{1-1/p,p}(\partial\Omega)$; then $u\in W^{2,p}(\Omega)$.
    }\\
    Let now $u$ be a Steklov eigenfunction; then $u$ is a solution of the problem above with $f=u$ and $g=\ld u$. We first apply the result for $p=2$ (note that since a Steklov eigenfunction is in $W^{1,2}(\Omega)$, its trace on $\partial\Omega$ is indeed in $W^{1/2,2}(\partial\Omega)$), obtaining that $u\in W^{2,2}(\Omega)$. Then, by Sobolev embedding, $u\in W^{1,p}(\Omega)$ for some $p>2$; it follows that $f\in L^p(\Omega)$ and $g\in W^{1-1/p,p}(\partial\Omega)$ for this $p$, and we apply again the theorem to get $u\in W^{2,p}(\Omega)$. It is now clear how to bootstrap the argument to obtain $u\in W^{2,p}(\Omega)$ for any $p<\infty$, and Sobolev embedding theorems then imply that $u\in\mathcal{C}^1(\overline{\Omega})$.
\end{proof}
\end{proposition}
While Proposition \ref{nodal} below is quite intuitive and mirrors what happens for Laplace eigenfunctions, the proof needs some care due to the fact that nodal domains may not be Lipschitz domains; a workaround to this was suggested in \cite{BM}, Appendix D, for the Dirichlet eigenfunction case, and then in \cite{LE} for Neumann and Robin eigenfunctions. We follow these two references in the proof below. 

\begin{remark}
     A word on notation: since a Steklov eigenfunction $u$ is $\mathcal{C}^1$ up to the boundary, in the following we implicitly consider an extension of it in a neighbourhood of $\overline{\Omega}$; when we talk about nodal domains of $u$, we are thinking about nodal domains of this extension. In this way the intersection of a nodal domain with $\partial\Omega$ is a relatively open set in $\partial\Omega$.
\end{remark}
\begin{proposition}
\label{nodal}
    Assume that $\partial\Omega$ is $\mathcal{C}^2$, Let $u=u_\ld$ be a solution of \eqref{problem}, and denote by $\mathcal{D}$ a nodal domain of $u$. Then $u$ is the first eigenfunction of a Steklov-Dirichlet problem in $\mathcal{D}$, with $S=\partial{\mathcal{D}}\cap\Omega$ and $F=\mathcal{D}\cap\partial\Omega$ in the notation of \eqref{mixed}; that is to say that 
    \begin{align}
    \label{nodalsteklov}
        \ld=\frac{\int_{\mathcal{D}}|\nabla u|^2}{\int_{F}|u|^2}=\inf \frac{\int_{\mathcal{D}}|\nabla \phi|^2}{\int_{F}|\phi|^2},
    \end{align}
    where the $\inf$ is taken over all $\phi\in\mathcal{C}^{\infty}(\overline{\mathcal{D}})$ which are compactly supported in $\mathcal{D}\cup F$.
    \begin{proof}
    As remarked earlier, every nodal domain of $u$ intersects $\partial \Omega$. Suppose that $u$ is positive in $\mathcal{D}$. If $\mathcal{D}$ were a Lipschitz domain, the proposition would immediately follow from the considerations in Section 2. Following \cite{BM} and \cite{LE} we consider instead the set $\mathcal{D}_{\alpha}=\{x\in\mathcal{D}: u(x)>\alpha\}$, for $\alpha$ small enough. If $\alpha$ is a regular value both of $u$ and of $u$ restricted to $\partial\Omega$ (recall that $u\in \mathcal{C}^1(\overline\Omega)$), then $\mathcal{D}_{\alpha}$ is a Lipschitz domain. Set $S_{\alpha}=\partial\mathcal{D}_{\alpha}\cap\Omega$, $F_{\alpha}=\mathcal{D}_{\alpha}\cap\partial\Omega$, and $u_{\alpha}=u-\alpha$. By Green's formula in $\mathcal{D_{\alpha}}$ and the Steklov condition we get that 
    \begin{align*}
        \int_{\mathcal{D}_{\alpha}}|\nabla u|^2=\int_{F_{\alpha}}u_{\alpha}\partial_{\nu}u_{\alpha}=\ld \int_{F_{\alpha}}|u|^2-\alpha\ld\int_{F_{\alpha}}u.
    \end{align*}
    Sard's Lemma gives a sequence of regular values for both $u$ and its restriction to the boundary approaching $0$, so that passing to the limit in the equality above we obtain 
    \begin{align}
    \label{green}
        \int_{\mathcal{D}}|\nabla u|^2=\ld \int_{F}|u|^2,
    \end{align}
    which gives the first equality in \eqref{nodalsteklov}. We consider now another auxiliary function in $\mathcal{D_{\alpha}}$, where $\alpha$ is a regular value as above; namely, we let $v_{\alpha}$ be the first eigenfunction of the mixed Steklov-Dirichlet problem in $\mathcal{D_{\alpha}}$, with Steklov condition on $F_{\alpha}$ and Dirichlet condition on $S_{\alpha}$. As in Section 2, we denote by $\sigma_1(\mathcal{D_{\alpha}})$ the corresponding first Steklov-Dirichlet eigenvalue. By Green's formula, we have that 
    \begin{align*}
        \int_{\mathcal{D_{\alpha}}}\nabla u \nabla v_{\alpha}=\int_{S_{\alpha}}u\partial_{\nu}v_{\alpha}+\int_{F_{\alpha}}u\partial_{\nu}v_{\alpha}=\int_{S_{\alpha}}u\partial_{\nu}v_{\alpha}+\sigma_1(\mathcal{D_{\alpha}})\int_{F_{\alpha}}uv_{\alpha}
    \end{align*}
    and also that 
    \begin{align*}
        \int_{\mathcal{D_{\alpha}}}\nabla u \nabla v_{\alpha}=\int_{S_{\alpha}}v_{\alpha}\partial_{\nu}u+\int_{F_{\alpha}}v_{\alpha}\partial_{\nu}u=\ld\int_{F_{\alpha}}uv_{\alpha},
    \end{align*}
    and therefore
    \begin{align*}
        \ld\int_{F_{\alpha}}uv_{\alpha}=\sigma_1(\mathcal{D_{\alpha}})\int_{F_{\alpha}}uv_{\alpha}+\int_{S_{\alpha}}u\partial_{\nu}v_{\alpha}.
    \end{align*}
    Without loss of generality, let $v_{\alpha}$ be positive in $\mathcal{D_{\alpha}}$; then on $S_{\alpha}$ we have $\partial_{\nu}v_{\alpha}\leq 0$, from which we obtain by the equation above that 
    \begin{align}
    \label{above}
        \sigma_1(\mathcal{D_{\alpha}})\geq \ld.
    \end{align}

Let us call $\sigma_1(\mathcal{D})$ the $\inf$ on the right hand side of \eqref{nodalsteklov}; since $\mathcal{D}$ may not be a Lipschitz domain, this is a slight abuse of notation. Using \eqref{green}, we have the string of inequalities
\begin{align*}
    \ld\int_{F}|u|^2=\int_{\mathcal{D}}|\nabla u|^2\geq \int_{\mathcal{D_{\alpha}}}|\nabla u_{\alpha}|^2\geq \sigma_1(\mathcal{D_{\alpha}})\int_{F_{\alpha}}|u_{\alpha}|^2\geq \sigma_1(\mathcal{D}) \int_{F_{\alpha}}|u_{\alpha}|^2,
\end{align*}
from which, letting $\alpha\to 0$, we obtain $\ld\geq \sigma_1(\mathcal{D})$. 

For the converse inequality, let $\varepsilon>0$ and pick $\phi$ compactly supported in $\mathcal{D}\cup F$ such that 
\begin{align*}
    \frac{\int_{\mathcal{D}}|\nabla \phi|^2}{\int_{F}|\phi|^2}\leq \sigma_1(\mathcal{D})+\varepsilon.
\end{align*}
We take now $\alpha$ a regular value as above so small that the support of $\phi$ is contained in $\mathcal{D}_{\alpha}\cup F_{\alpha}$; note that here we can exclude the case of a point $x_0\in\partial\Omega$ such that $u(x_0)=0$ and $u>0$ in a neighbourhood of $x_0$ by Hopf's lemma combined with the Steklov condition. Using \eqref{above}, we then have 
\begin{align*}
    \ld\leq \sigma_1(\mathcal{D_{\alpha}})\leq \frac{\int_{\mathcal{D}}|\nabla \phi|^2}{\int_{F}|\phi|^2}\leq \sigma_1(\mathcal{D})+\varepsilon.
\end{align*}
Since this holds for any $\varepsilon>0$, we finally obtain $\ld=\sigma_1(\mathcal{D})$, which gives the second equality in \eqref{nodalsteklov} and concludes the proof of the proposition. 
\end{proof}
\end{proposition}

The proof of Theorem \ref{upper} is now quite easy, and the main ingredient is once again provided by Lemma \ref{key}. Suppose that a nodal domain of $u$, $\mathcal{D}$, contains a set $W=B(x_0,s)\cap \Omega$, for some $x_0\in \partial \Omega$ and some $s>0$. We assume $s$ is small enough to apply Lemma \ref{key}. Consider problem \eqref{mixed} on $W$, with $S=\partial B(x_0,s)\cap \mathcal{D}$ and $F=\partial\Omega\cap B(x_0,s)$; by Lemma \ref{key}, we have $\sigma_1(W)<C/s$. On the other hand, by Proposition \ref{nodal}, $\ld=\sigma_1(\mathcal{D})$, where as in the proof of Proposition \ref{nodal} by $\sigma_1(\mathcal{D})$ we denote the $\inf$ on the right hand side of \eqref{nodalsteklov}. This definition immediately implies weak domain monotonicity, $\sigma_1(W)\geq \sigma_1(\mathcal{D})$. We obtain
\begin{align*}
    C/s>\sigma_1(W)\geq \sigma_1(\mathcal{D})=\ld,
\end{align*}
so that $s<C/\ld$, and the proof is once again complete.

\section{Lower bound on the inner radius}
We discuss here the proof of Theorem \ref{lower}. The argument is quite elementary, and the ideas follow those for Laplace eigenfunctions in \cite{M} and especially \cite{M2}; the boundary presents some additional difficulties. 

From now on, $\Omega$ will be a Lipschitz domain in $\R^2$. We start by covering $\partial\Omega$ with a finite number $K=K(\Omega)$ of (possibly overlapping) rectangles centered at $\partial\Omega$, such that in each such rectangle $\Omega$ is the set of points below the graph of a Lipschitz function; let us denote the Lipschitz constant of this function by $\Gamma$. In the sequel we will refer to such rectangles as "large rectangles". After an orthogonal transformation we have the following description of $\Omega$ in a large rectangle $Q$:
\begin{align}
\label{lip}
    \Omega\cap Q=\{(x,y)\in Q:-1<x<1 , -2\Gamma< y < g(x)\},
\end{align}
where $g$ is a Lipschitz function with $g(0)=0$, $|g(x_1)-g(x_2)|\leq \Gamma |x_1-x_2|$ for $x_1, x_2 \in [-1,1]$. The boundary is described by $\partial\Omega\cap Q=\{(x,y)\in Q: y=g(x)\}$. Note that we choose the rectangles in a way that the graph of $g$ intersects the boundary of each rectangle only on the vertical sides.\\

Let now $u_{\ld}$ be a Steklov eigenfunction and $\mathcal{D}$ any of its nodal domains. As in Section 3, we use the notation $F=\partial\Omega\cap \overline{\mathcal{D}}$ and $S=\mathcal{Z}_{u_{\ld}}\cap\overline{\mathcal{D}}$ (note that by our definition $\mathcal{Z}_{u_{\ld}}\subset\Omega$). We have that 
\begin{align}
\label{ray}
    \ld=\frac{\int_{\mathcal{D}}|\nabla u_{\ld}|^2}{\int_F |u_{\ld}|^2}.
\end{align}

As in the Introduction, we denote by $\rho(\mathcal{D})$ the modified inner radius of $\mathcal{D}$; we repeat that by this we mean the radius of the largest ball $B$ centered at $F$ such that $B\cap\Omega\subset \mathcal{D}$. Fix 
\begin{align*}
    h=100\rho(\mathcal{D}); 
\end{align*}
we can assume that $h$ is small enough. For each rectangle $Q$ in the above covering of $\partial\Omega$ we cover $F\cap Q$ by non-overlapping almost rectangles $q$; by almost rectangle we mean a set of the form 
\begin{align}
\label{cubes}
    q=\{(x,y)\in Q:x_q-h<x<x_q+h, g(x)-2\Gamma h < y < g(x)+2\Gamma h\},
\end{align}
for some $x_q\in (-1+h,1-h)$ and $g$ as above. Again the choice of the sides of the almost rectangle is made to ensure that the graph of $g$ intersects each $q$ on the vertical sides only.
Let now 
\begin{align*}
    u=
    \begin{cases}
      u_{\ld} \quad \text{on } \mathcal{D};\\
      0 \quad \ \; \text{on } \Omega\setminus\mathcal{D}.
    \end{cases}
\end{align*}    
Using \eqref{ray} and remembering that almost rectangles do not overlap, it is easy to verify the following fact.
\begin{claim}
\label{special}
There exists at least one almost rectangle $q$ such that 
\begin{align*}
    \int_{q\cap\Omega}|\nabla u|^2\leq K\ld \int_{q\cap F}|u|^2,
\end{align*}
where $K$ is as above the number of large rectangles in the covering of $\partial\Omega$.
\end{claim}

If we can show that in the almost rectangle given in Claim \ref{special} we have 
\begin{align}
\label{poin}
    \int_{q\cap F}|u|^2\leq Ch \int_{q\cap\Omega}|\nabla u|^2,
\end{align}
with constant $C$ depending only on $\Omega$ (actually on $\Gamma$), Theorem \ref{lower} follows. From now on, the letter $C$ will denote constants that may change value at every occurrence but depend on $\Gamma$ only. We will assign explicit numerical values to constants independent of $\Gamma$. 

We now prove \eqref{poin} using some Poincaré-type inequalities, in the spirit of \cite{M2}. Let us first note a trace inequality; the proof could be obtained by general Sobolev trace inequalities and scaling, but since it is simple we give it here for the reader's convenience. We denote 
\begin{align}
\label{map}
    \varphi(x,y)=y-g(x).
\end{align}
In the following, for brevity, we consider $\varphi$ as defined in $q$, so that we denote by $\{\varphi=t\}$ the set $\{(x,y)\in q : \varphi(x,y)=t\}$.
\begin{lemma}
In the notation above, we have that 
\begin{align}
\label{trace}
    \int_{q\cap F}|u|^2\leq C\left(\frac{1}{h}\int_{q\cap\Omega}|u|^2+h\int_{q\cap\Omega}|\nabla u|^2\right),
\end{align}
where $C$ is independent of $h$.
\begin{proof}
By the co-area formula and the Lipschitz bound of $g$,
\begin{align*}
    \int_{q\cap\Omega}|u|^2\geq \frac{1}{1+\Gamma}\int_{q\cap\Omega}|u|^2|\nabla \varphi|=\frac{1}{1+\Gamma}\int_{-2\Gamma h}^{0}\left(\int_{\{\varphi=t\}}|u|^2\right)dt,
\end{align*}
from which it follows that there is a $t_0\in(-2\Gamma h,0)$ such that 
\begin{align*}
    \int_{\{\varphi=t_0\}}|u|^2\leq \frac{C}{h} \int_{q\cap\Omega}|u|^2,
\end{align*}
with $C$ depending on $\Gamma$. Take now a point $(x,g(x))\in q\cap F$; we have that
\begin{align*}
    |u(x,g(x))|^2\leq 2|u(x,g(x)+t_0)|^2+2|u(x,g(x))-u(x,g(x)+t_0)|^2
    \\
    \leq 2|u(x,g(x)+t_0)|^2+Ch\int_{-2\Gamma h}^{0}|\nabla u|^2(x,g(x)+t)dt.
\end{align*}
Integrating along $q\cap F$, we obtain
\begin{align*}
    \int_{q\cap F}|u|^2\leq 2 \int_{\{\varphi=t_0\}}|u|^2+ Ch\int_{q\cap\Omega}|\nabla u|^2 \leq C\left(\frac{1}{h}\int_{q\cap\Omega}|u|^2+h\int_{q\cap\Omega}|\nabla u|^2\right).
\end{align*}
\end{proof}
\end{lemma}

The almost rectangle $q$ identified in Claim \ref{special} will be fixed; we call the point $(x_q,g(x_q))$ the center of $q$. We denote by $(1/4)q$ the almost rectangle with center $(x_q,g(x_q))$ and with side lengths the ones of $q$ divided by $4$. By the choice of $h$, we have that 
\begin{align}
\label{deep}
    \frac{1}{4}q\cap (\Omega\setminus\mathcal{D})\neq \emptyset.
\end{align}
(Otherwise, $(1/4)q\cap\Omega \subset \mathcal{D}$ and hence $\mathcal{D}$ would contain the intersection with $\Omega$ of a ball of radius $2\rho(\mathcal{D})$, which contradicts the definition of modified inner radius). Take now any point $p$ in the above intersection; following the notation of \cite{M2}, we denote by $hole(p)$ the nodal domain of $u_{\ld}$ that contains $p$. By the maximum principle, $hole(p)$ must intersect $\partial(q\cap\Omega)$. There are essentially 3 configurations that can arise; we will analyse them separately.

\subsection*{Configuration 1:} $hole(p)\subset q\cap \overline{\Omega}$. This means that $hole(p)$ exits $q\cap\Omega$ through $\partial\Omega$, that is it intersects $\partial(q\cap\Omega)$ only on $\partial\Omega$. As in Section 2 and 3, we denote by $\sigma_1(q\cap\Omega)$ the first eigenvalue of the Steklov-Dirichlet problem with Steklov condition on $q\cap F$, Dirichlet on $\partial(q\cap\Omega)\setminus F$. The following claim is proved with an easy Poincaré inequality.
\begin{claim}
The estimate 
\begin{align*}
    \sigma_1(q\cap\Omega)\geq \frac{c}{h}
\end{align*}
holds with a constant $c>0$ independent of $h$.
\begin{proof}
Let $v$ be any function in $H^1(q\cap\Omega)$ which has zero trace on $\partial(q\cap\Omega)\setminus F$. In particular, with $\varphi$ as in \eqref{map}, $v=0$ on $\{\varphi=-2\Gamma h\}$. For $(x,g(x))\in F$, we then have 
\begin{align*}
    |v(x,g(x))|^2=|v(x,g(x))-v(x,g(x)-2\Gamma h)|^2=|\int_{-2\Gamma h}^{0} \partial_y v(x,g(x)+t)dt|^2\\
    \leq 2\Gamma h\int_{-2\Gamma h}^{0}|\nabla v|^2 (x,g(x)+t)dt.
\end{align*}
Integrating along $F$, we obtain 
\begin{align*}
    \int_{q\cap F}|v|^2\leq Ch \int_{q\cap\Omega}|\nabla v|^2;
\end{align*}
the claim now follows by recalling the variational characterization of $\sigma_1$ in \eqref{eigenvalue}.
\end{proof}
\end{claim}
Since $hole(p)$ is a nodal domain of $u_{\ld}$, we have that
\begin{align*}
    \ld=\frac{\int_{hole(p)}|\nabla u_{\ld}|^2}{\int_{\overline{hole(p)}\cap F}|u_{\ld}|^2}.
\end{align*}
Recall that $hole(p)\subset q\cap \overline{\Omega}$; the function $u_{\ld}\chi_{hole(p)}$ is therefore in competition in the inf defining $\sigma_1(q\cap\Omega)$. By the above claim, we then obtain 
\begin{align*}
    h\geq \frac{c}{\ld}.
\end{align*}
Recalling that $h=100\rho(\mathcal{D})$, Theorem \ref{lower} is proved in this configuration. 

\subsection*{Configuration 2:} $hole(p)$ intersects $\partial(q\cap\Omega)$ only in $\{\varphi=-2\Gamma h\}$. In this situation, since $p\in (1/4)q$, every curve of the form $\{\varphi=t\}$ for $t\in(-2\Gamma h,-\Gamma h)$ intersects the zero set of $u$. Consider the set 
\begin{align*}
    E_1=\{(x,y)\in q: g(x)-2\Gamma h\leq y \leq g(x)-\Gamma h\}.
\end{align*}
\begin{claim}
Under the above hypothesis,
\begin{align*}
    \int_{E_1}|u|^2\leq Ch^2 \int_{E_1}|\nabla u|^2.
\end{align*}
\begin{proof}
    By the co-area formula,
    \begin{align*}
        \int_{E_1}|u|^2\leq \int_{-2\Gamma h}^{-\Gamma h}\left(\int_{\{\varphi=t\}}|u|^2\right)dt;
    \end{align*}
    since for any $t\in(-2\Gamma h,-\Gamma h)$ there is a point $p_t$ on the curve $\{\varphi=t\}$ such that $u(p_t)=0$, a change of variables and a 1-dimensional Poincaré inequality give
    \begin{align*}
        \int_{\{\varphi=t\}}|u|^2\leq Ch^2 \int_{\{\varphi=t\}}|\nabla u|^2.
    \end{align*}
    Plugging it in the above and using again the co-area formula we obtain 
    \begin{align*}
        \int_{E_1}|u|^2\leq Ch^2\int_{-2\Gamma h}^{-\Gamma h}\left(\int_{\{\varphi=t\}}|\nabla u|^2\right)dt\leq Ch^2\int_{E_1}|\nabla u|^2|\nabla \varphi|\leq C'h^2 \int_{E_1}|\nabla u|^2.
    \end{align*}
\end{proof}
\end{claim}

Note that there is at least one $t_0\in(-2\Gamma h,-\Gamma h)$ such that 
    \begin{align*}
        \int_{\{\varphi=t_0\}}|u|^2\leq \frac{C}{h}\int_{E_1}|u|^2.
    \end{align*}
We now procede as in the proof of \eqref{trace}: for $(x,g(x))\in q\cap F$, 
\begin{align*}
    |u(x,g(x))|^2\leq 2|u(x,g(x)+t_0)|^2+2|u(x,g(x))-u(x,g(x)+t_0)|^2
    \\
    \leq 2|u(x,g(x)+t_0)|^2+Ch\int_{-2\Gamma h}^{0}|\nabla u|^2(x,g(x)+t)dt.
\end{align*}
Integrating along $F$ and using Claim 3 we obtain 
\begin{align*}
    \int_{q\cap F}|u|^2\leq 2\int_{\{\varphi=t_0\}}|u|^2+Ch\int_{q\cap\Omega}|\nabla u|^2\leq \frac{C}{h}\int_{E_1}|u|^2+Ch\int_{q\cap\Omega}|\nabla u|^2\\
    \leq Ch\int_{E_1}|\nabla u|^2+Ch\int_{q\cap\Omega}|\nabla u|^2\leq Ch \int_{q\cap\Omega}|\nabla u|^2,
\end{align*}
which is the inequality \eqref{poin} we needed to show to prove Theorem \ref{lower}.

\subsection*{Configuration 3:} $hole(p)$ intersects $\partial(q\cap\Omega)$ in at least one of the vertical sides of $q$; say it crosses the left side. Since $p\in(1/4)q$, we have that every line segment of the form \begin{align*}
    I_s=\{(x,y)\in q\cap \Omega: x=s\}
\end{align*}
for $s\in (x_q-h,x_q-h/2)$ intersects the zero set of $u$. Let
\begin{align*}
    E_2=\bigcup_{s\in (x_q-h,x_q-h/2)}I_s.
\end{align*} 
The following claim is proved in the same way as Claim 3.
\begin{claim}
\begin{align*}
    \int_{E_2}|u|^2\leq Ch^2\int_{E_2}|\nabla u|^2.
\end{align*}
\end{claim}
By the co-area formula we have that there exists some $s_0$ such that \begin{align*}
    \int_{I_{s_0}}|u|^2\leq \frac{C}{h}\int_{E_2}|u|^2.
\end{align*}
Consider now any point $p\in q\cap\Omega$, which we write as $p=(s,g(s)-t)$; let $p_0=(s_0,g(s_0)-t)\in I_{s_0}$. Note that both $p$ and $p_0$ lie on the curve $\{\varphi=t\}$. We have that 
\begin{align*}
    |u(p)|^2\leq 2|u(p_0)|^2+2|u(p)-u(p_0)|^2\leq 2|u(p_0)|^2+Ch\int_{\{\varphi=t\}}|\nabla u|^2.
\end{align*}
We integrate the above over $q\cap\Omega$ using the co-area formula:
\begin{align*}
    \int_{q\cap\Omega}|u(p)|^2\leq \int_{-2\Gamma h}^0 \left(\int_{\{\varphi=t\}}|u(p)|^2\right)dt\leq Ch\int_{-2\Gamma h}^0|u(p_0)|^2dt \\ + Ch^2\int_{-2\Gamma h}^0 \left(\int_{\{\varphi=t\}}|\nabla u|^2\right)dt.
\end{align*}
By the choice of $s_0$ and Claim 4, we then obtain 
\begin{align*}
    \int_{q\cap\Omega}|u|^2\leq Ch^2\int_{q\cap\Omega}|\nabla u|^2.
\end{align*}
Finally, by the above and the trace inequality \eqref{trace}, we get 
\begin{align*}
    \int_{q\cap F}|u|^2\leq Ch \int_{q\cap\Omega}|\nabla u|^2.
\end{align*}
This is once again the desired inequality \eqref{poin}, which concludes our analysis. Theorem \ref{lower} is proved. 

\subsection*{Comments about higher dimensions} The proof works in dimension two because the length of some projection of $hole(p)$ on a side of the almost rectangle is large enough (independent of $\ld$). In higher dimensions we cannot quite say the same. For Laplace eigenfunctions, a partial workaround is used in \cite{M}: one estimates the local asymmetry of an eigenfunction, that is the volume of the set of positivity in a ball, and then uses a Poincaré inequality of Maz'ya \cite{MZ} involving capacities. This gives a lower bound on the inner radius for Laplace eigenfunctions that is polynomial in the reciprocal of the eigenvalue (see \cite{M} for the dimensional exponent); as far as we know, in dimension greater than two the conjectured optimal lower bound has not been obtained. 
A similar approach can be used for Steklov eigenfunctions, if the boundary is regular enough; in this case one can extend a Steklov eigenfunction across the boundary, which facilitates the proof of a local asymmetry estimate for the eigenfunction near the boundary. One can then prove a Poincaré inequality in the spirit of \cite{MZ} and obtain a lower bound on the inner radius that is polynomial in $\ld^{-1}$ (with an explicit dimensional exponent). As mentioned in the introduction, we choose not to give the details in the present work and reserve it to subsequent work. 

\section*{Acknowledgements}
This work owes a lot to the patient guidance of Eugenia Malinnikova; discussing the article with her has contributed greatly to both the contents and the presentation, and her reading of various drafts helped spot mistakes and inaccuracies. Many thanks are due to Aleksandr Logunov for his encouragement, as well as for reading a previous draft and asking useful questions. Iosif Polterovich read this work at a previous stage and provided a wealth of helpful comments, for which I am very grateful. I also wish to thank the reviewers for their thorough reading and several suggestions that helped improve the presentation.

This work was started while I was a Visiting Student Researcher at the Department of Mathematics at Stanford University; it is a pleasure to thank the department for the hospitality and the nice working conditions.

The author is supported by Project 275113 of the Research Council of Norway.


\end{document}